\documentclass[reqno]{amsart}

\pdfoutput = 1

\usepackage[english]{babel}
\usepackage{amsmath,amssymb}
\usepackage{amsfonts}
\usepackage{array,dcolumn}
\usepackage{graphicx}

\usepackage[T1]{fontenc}
\usepackage[utf8]{inputenc}
\usepackage[activate={true,nocompatibility},spacing,kerning]{microtype}
\usepackage[hyperpageref]{backref}
\usepackage[pdftex,
            colorlinks,
            hyperfootnotes=false,
            pdffitwindow=true,
            plainpages=false,
            pdfpagelabels=true,
            pdfpagemode=UseOutlines,
            pdfpagelayout=SinglePage,
            pdftitle={Optimal Contours for High-Order Derivatives},
            pdfauthor={Folkmar Bornemann, Georg Wechslberger, Technische Universität München},
            pagebackref,
            hyperindex,
            filecolor=blue,
            urlcolor=blue,
            citecolor=blue,
            linkcolor=blue]{hyperref}
\usepackage{myharvard}
\usepackage{bm}
\usepackage[sc,osf]{mathpazo}
\graphicspath{{plots/}}

\renewcommand{\leq}{\leqslant}
\renewcommand{\geq}{\geqslant}

\newcommand{\SP}{\mathcal{P}}
\let\join\relax
\DeclareMathOperator*{\argmin}{argmin}
\DeclareMathOperator{\join}{\dot{+}}
\DeclareMathOperator{\ind}{ind}
\DeclareMathOperator{\Ai}{Ai}

\begin{document}

\renewcommand{\figurename}{\footnotesize {\sc Figure\@}}
\renewcommand{\tablename}{\footnotesize {\sc Table\@}}

\title{Optimal Contours for High-Order Derivatives}

\author{Folkmar Bornemann}
\author{Georg Wechslberger}
\address{Zentrum Mathematik -- M3, Technische Universität München,
         80290~München, Germany}
\email{bornemann@tum.de; wechslbe@ma.tum.de}


\begin{abstract} As a model of more general contour integration problems
we consider the numerical calculation of high-order derivatives of holomorphic functions using Cauchy's integral formula.
\citeasnoun{springerlink:10.1007/s10208-010-9075-z} showed that the condition number of the Cauchy
integral strongly depends on the chosen contour and solved the
problem of minimizing the condition number for circular contours. In this paper we
minimize the condition number within the class of grid paths of step size $h$ using
Provan's algorithm for finding a \emph{shortest enclosing walk} in weighted graphs embedded in the plane.
Numerical examples show that optimal grid paths yield small condition numbers
even in those cases where circular contours are known to be of limited use, such as for functions with
branch-cut singularities.
\end{abstract}

\maketitle

\section{Introduction}\label{sec:intro}

To escape from the ill-conditioning of difference schemes for the numerical
calculation of high-order derivatives, numerical quadrature applied to
Cauchy's integral formula has on various occasions been suggested as a remedy
\citeaffixed{springerlink:10.1007/s10208-010-9075-z}{for a survey of the literature, see}.
To be specific, we consider a function $f$ that is holomorphic on a complex domain $D \ni 0$; Cauchy's formula
gives\footnote{Without loss of generality we evaluate derivatives at $z=0$.}
\begin{equation}\label{eq:cauchy}
f^{(n)}(0) = \frac{n!}{2\pi i} \int_{\Gamma} z^{-n-1} f(z)\,dz
\end{equation}
for each cycle $\Gamma \subset D$ that has winding number $\ind(\Gamma;0)=1$. If $\Gamma$ is not carefully chosen,
however, the integrand tends to oscillate at a frequency of order $O(n^{-1})$ with very large amplitude
\cite[Fig.~4]{springerlink:10.1007/s10208-010-9075-z}. Hence, in general, there is much cancelation in the
evaluation of the integral and ill-conditioning returns through the backdoor. The condition number of the integral\footnote{Given an accurate and stable (i.e., with positive weights)
quadrature method such as Gauss--Legendre or Clenshaw--Curtis, this condition number actually yields, by 
\[
\text{\# loss of significant digits} \approx \log_{10} \kappa(\Gamma,n),
\]
an estimate of the error caused by round-off in the last significant digit of the data (i.e., the function $f$).}
is \cite[Lemma~9.1]{MR1949263}
\[
\kappa(\Gamma,n) = \frac{\int_{\Gamma} |z|^{-n-1} |f(z)|\,d|z|}{\left|\int_{\Gamma} z^{-n-1} f(z)\,dz\right|}
\]
and $\Gamma$ should be chosen as to make this number as small as possible. Equivalently, since the denominator is, by
Cauchy's theorem, independent of $\Gamma$, we have to minimize
\begin{equation}\label{eq:weight}
d(\Gamma) = \int_{\Gamma} |z|^{-n-1} |f(z)|\,d|z|.
\end{equation}
\citeasnoun{springerlink:10.1007/s10208-010-9075-z} considered circular contours of radius $r$; he found that there is a unique $r_*=r(n)$ solving
the minimization problem and that there are different scenarios for the corresponding condition number $\kappa_*(n)$ as $n\to \infty$:
\begin{itemize}
\item $\kappa_*(n) \to \infty$, if $f$ is in the Hardy space $H^1$;\\*[-3.5mm]
\item $\limsup_{n\to \infty} \kappa_*(n) \leq M$, if $f$ is an entire function of completely regular growth
which satisfies a non-resonance condition of
the zeros and whose Phragmén--Lindelöf indicator possesses $M$ maxima (a small integer).
\end{itemize}
Hence, though those (and similar) results basically solve the problem of choosing proper contours for entire functions,
much better contours have to be found for the class $H^1$. Moreover, the restriction to circles lacks any algorithmic flavor that would point to more general problems depending on the
choice of contours, such as the numerical solution of highly-oscillatory Riemann--Hilbert problems \cite{Olver2011}.\footnote{Taking the contour optimization developed in this paper as a model, \citeasnoun{1206.2446} has recently
addressed the deformation of Riemann--Hilbert problems from an algorithmic
point of view.} 

\begin{figure}[tbp]
\begin{center}
\includegraphics[width=0.425\textwidth]{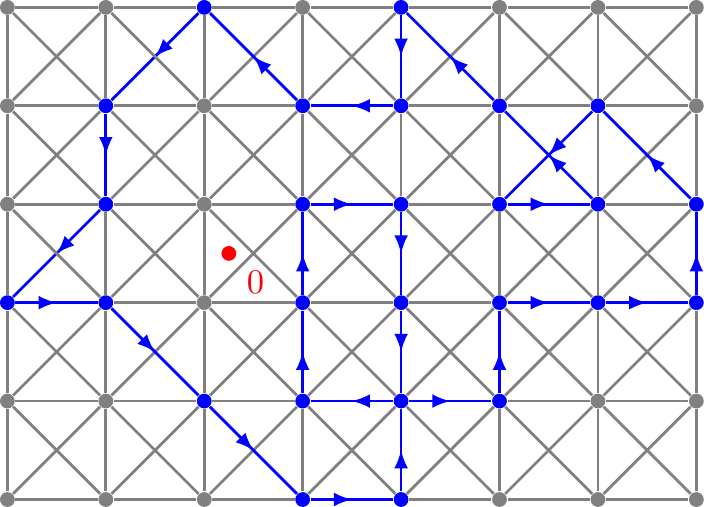}
\end{center}
\caption{\footnotesize Path $\Gamma$ with $\ind(\Gamma;0)=1$ in a grid-graph of step size $h$.}\label{fig:rect}
\end{figure}

In this paper, we solve the contour optimization problem within the more general class of grid paths of step size $h$ (see Fig.~\ref{fig:rect}; we allow diagonals to be included)
as they are known from Artin's proof of the general, homological version of Cauchy's integral theorem \cite[IV.3]{MR1659317}. Such paths are composed from horizontal, vertical and diagonal edges taken from a (bounded) grid $\Omega_h \subset D$ of step size $h$. Now, the weight function (\ref{eq:weight}), being \emph{additive} on the abelian group of path chains, turns
the grid $\Omega_h$ into an edge-weighted graph such that each optimal grid path $W_*$ becomes a \emph{shortest enclosing walk} (SEW); ``enclosing''
because we have to match the winding number condition $\ind(W_*;0)=1$. An efficient solution of the SEW problem for embedded graphs was found by \citeasnoun{Provan1989119} and serves as a starting point for our work.

\subsection*{Outline of the Paper}

In Section~\ref{sec:graph} we discuss general embedded graphs in which an optimal contour is to be searched for; we discuss the problem of finding a
shortest enclosing walk and recall Provan's algorithm. In Section~\ref{sec:detail} we discuss some implementation details and tweaks for the problem at hand.
Finally, in Section~\ref{sec:numer} we give some numerical examples; these can easily be constructed
in a way that the new algorithm outperforms, by orders of magnitude, the optimal circles of \citeasnoun{springerlink:10.1007/s10208-010-9075-z} with respect to accuracy and the direct symbolic differentiation with respect to efficiency. 

\section{Contour Graphs and Shortest Enclosing Walks}\label{sec:graph}

By generalizing the grid $\Omega_h$, we consider a finite graph $G = (V,E)$ \emph{embedded to}~$D$, that is, built from vertices $V \subset D$ and edges $E$ that are smooth curves connecting the
vertices within the domain $D$. We write $uv$ for the edge connecting the vertices $u$ and $v$; by (\ref{eq:weight}), its weight is defined as
\begin{equation}\label{eq:edgeWeight}
d(uv)=\int_{uv} |z|^{-n-1} |f(z)|\,d|z|.
\end{equation}
A \emph{walk} $W$ in the graph $G$ is a \emph{closed} path built from a sequence of adjacent edges,
written as (where $\join$ denotes joining of paths)
\[
W= v_1v_2 \join v_2v_3 \join \cdots \join v_m v_1;
\]
it is called {\em enclosing} the obstacle $0$ if
the winding number is $\ind(W;0)=1$. The set of all possible enclosing walks is denoted by $\Pi$. As discussed in §\ref{sec:intro}, the condition number
is optimized by the shortest enclosing walk (not necessarily unique)
\[
W_* = \argmin_{W \in \Pi} d(W)
\]
where, with $W = v_1v_2 \join v_2v_3 \join \cdots \join v_m v_1$ and $v_{m+1}=v_1$, the \emph{total weight} is
\[
d(W) = \sum_{j=1}^m d(v_jv_{j+1}).
\]
The problem of finding such a SEW
was solved by \citeasnoun{Provan1989119}: the idea is that with $\SP_{u,v}$ denoting a shortest path between $u$ and $v$, any shortest enclosing walk
$W_* = w_1w_2 \join w_2w_3 \join \cdots \join w_m w_1$
can be cast in the form \cite[Thm.~1]{Provan1989119}
\begin{equation*}
	W_* = \SP_{w_1,w_j} \join w_jw_{j+1} \join \SP_{w_{j+1},w_1}
\end{equation*}
for at least one $j$. Hence, any shortest enclosing walk $W_*$ is already specified by one of its vertices and one of its edges; therefore
\begin{equation*}
	W_* \in \tilde\Pi = \{ \SP_{u,v} \join vw \join \SP_{w,u} : u \in V, vw \in E \}.
\end{equation*}
Provan's algorithm finds $W_*$ by, first, building the finite set $\tilde\Pi$; second, by removing all walks from it that do not enclose $z=0$; and third, by selecting a walk
from the remaining candidates that has the lowest total weight. Using \possessivecite{MR904195} implementation of Dijkstra's algorithm to compute the shortest paths $\SP_{u,v}$,
the run time of the algorithm is known to be \cite[Corollary~2]{Provan1989119}
\begin{equation}\label{eq:complexity}
O(|V|\,|E|+|V|^2\log|V|).
\end{equation}


\section{Implementation Details}
\label{sec:detail}

We restrict ourselves to graphs $\Omega_h$ given by finite square grids of step size $h$, centered at $z=0$---with all vertices and
edges removed that do not belong to the domain $D$. Since Provan's algorithm just requires an embedded graph but not a planar
graph, we may add the diagonals of the grid cells as further edges to the graph (see Fig.~\ref{fig:rect}).\footnote{These diagonals increase the number
of possible slopes which results, e.g., in improved approximations of the direction of steepest descent at a saddle point of $d(z)$ \cite[§9]{springerlink:10.1007/s10208-010-9075-z} or in a faster U-turn around the end of a branch-cut, see Fig.~\ref{fig:SewVsCircle}. The latter case leads to some significant reductions of the condition number, 
see Fig.~\ref{fig:nodes}.} 
For such a graph $\Omega_h$, with or without diagonals, we have $|V| = O(h^{-2})$ and $|E| = O(h^{-2})$ so that the
complexity bound (\ref{eq:complexity}) simplifies to $$O(h^{-4} \log h^{-1}).$$

\subsection{Edge Weight Calculation}

Using the edge weights $d(uv)$ on $\Omega_h$ requires to approximate the integral in \eqref{eq:edgeWeight}.
Since not much accuracy is needed here,\footnote{Recall that optimizing the condition number is just a question of order of magnitude but not of precise numbers. Once the contour $\Gamma$ has been fixed, a much more accurate quadrature rule will be employed to
calculate the integral (\ref{eq:cauchy}) itself, see §\ref{subsect:quadrature}.} a simple trapezoidal rule with two nodes is generally sufficient:
\begin{align*}
d(uv) &= \int_{uv} |z|^{-(n+1)} |f(z)| d |z| \\*[1mm]
	  &= \frac{|u - v|}{2} \left( d(u) + d(v) \right) + O(h^3) = \tilde d(uv) + O(h^3)
\end{align*}
with the \emph{vertex weight}
\begin{equation}\label{eq:vertexweight}
d(z) = |z|^{-(n+1)} |f(z)|.
\end{equation}
Although $\tilde d(uv)$ will typically have an accuracy of not more than just a few bits for the rather coarse
grids $\Omega_h$ we work with, we have not encountered a single case in which
a more accurate computation of the weights would have resulted in a different SEW $W_*$.

\begin{figure}[tbp]
\begin{center}
\begin{minipage}{0.35\textwidth}
\begin{center}
\includegraphics[width=\textwidth]{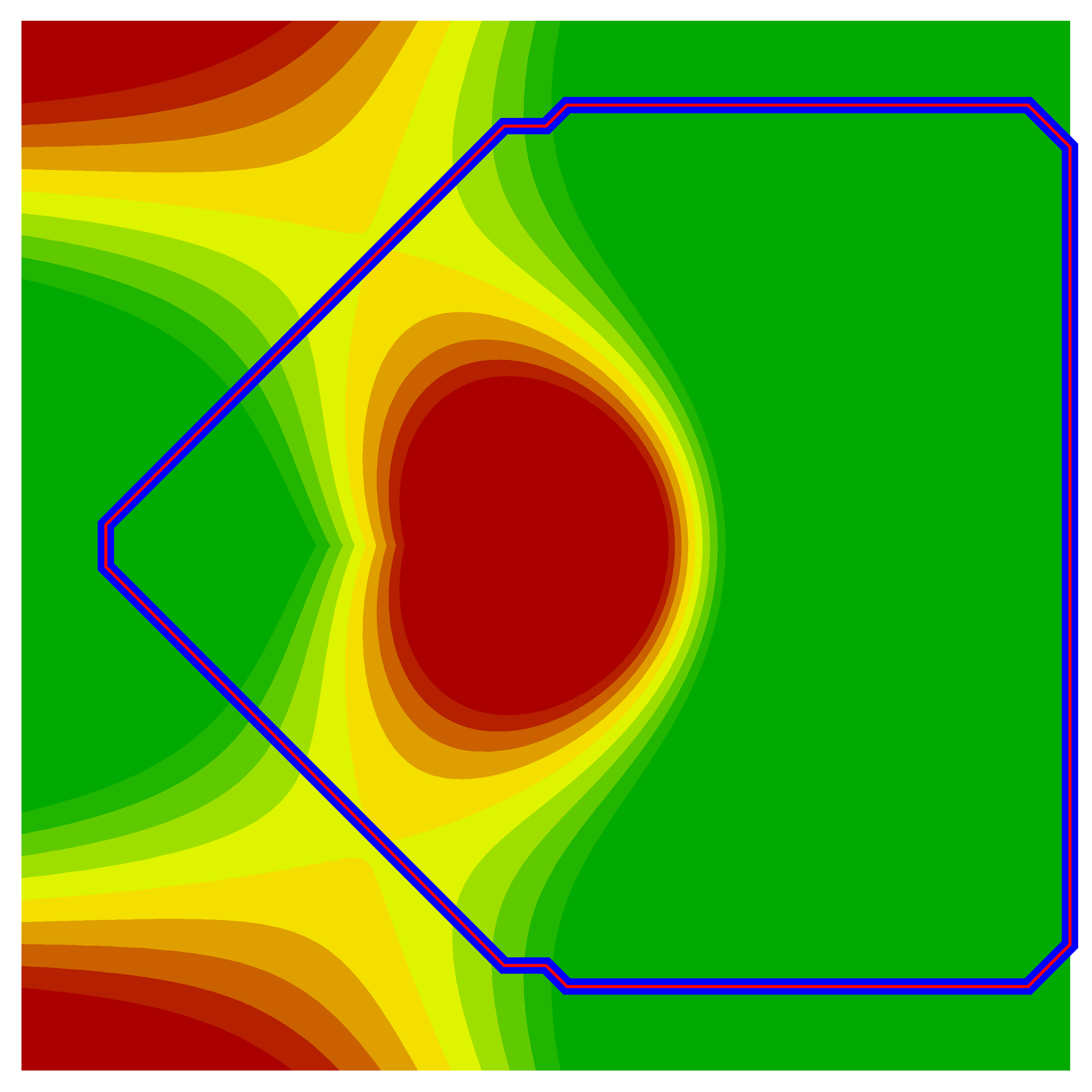}\\*[0.5mm]
{\footnotesize $\Ai(z)$}
\end{center}
\end{minipage}
\hfil
\begin{minipage}{0.35\textwidth}
\begin{center}
\includegraphics[width=\textwidth]{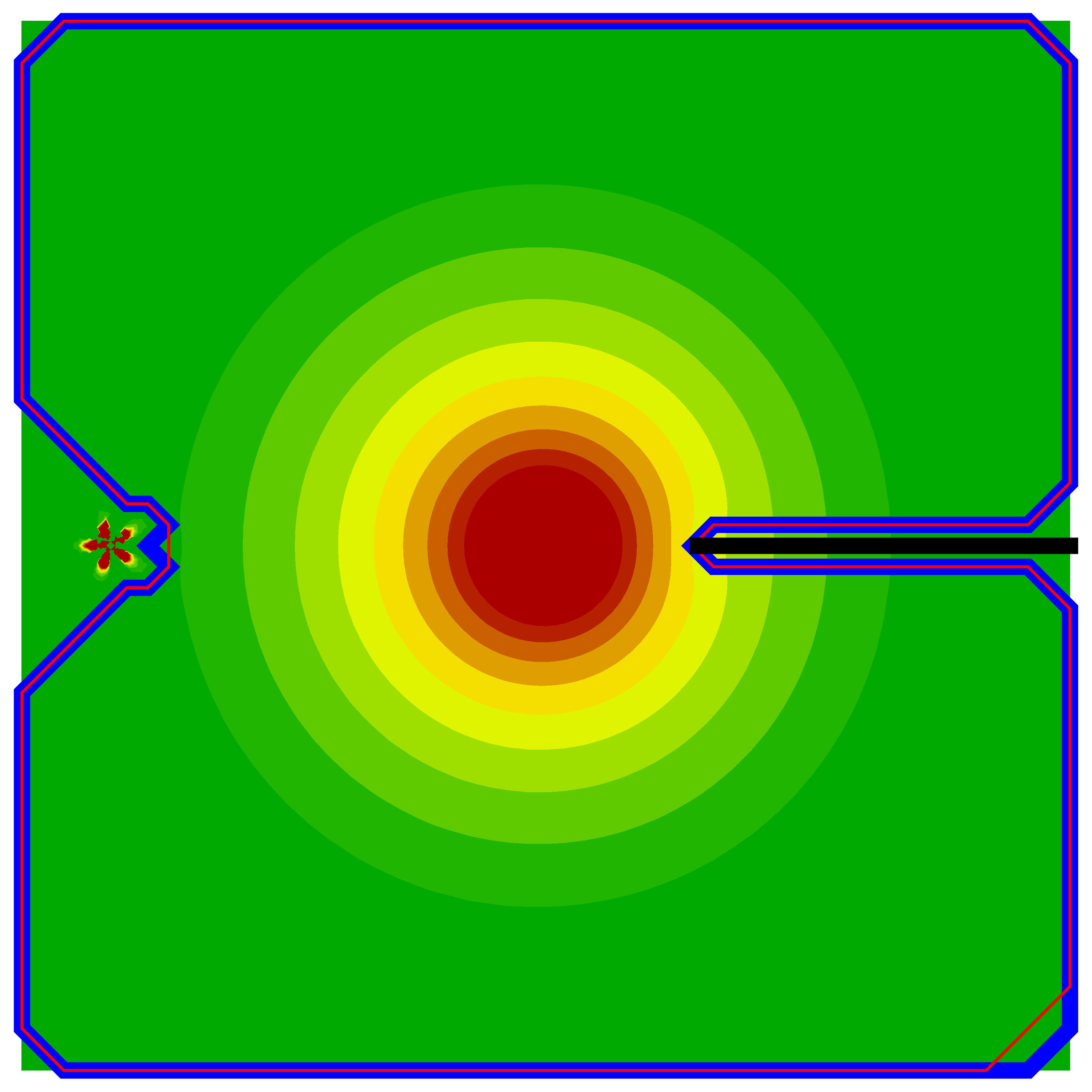}\\*[0.5mm]
{\footnotesize $\exp(1/(1+8z)^{1/5})(1-z)^{11/2}J_0(z)$}
\end{center}
\end{minipage}
\end{center}
	\caption{\footnotesize $W_*$ (red) vs. $W_{v_*}$ (blue): the color coding shows the size of $\log d(z)$; with red for large values and green for small values. The smallest level shown is the threshold, below of which the edges of $W_*$ do not contribute to the first couple of
significant digits of the total weight. The plots illustrate that $W_*$ and  $W_{v_*}$ differ typically just in a small region well below this threshold; consequently, both walks yield about
the same condition number. On the right note the five-leaved clover that represents the combination of algebraic
and essential singularity at $z=-1$.}
\label{fig:OneVsAll}
\end{figure}

\subsection{Reducing the size of $\tilde \Pi$}
As described in Section~\ref{sec:graph}, Provan's algorithm starts by building a walk for every pair $(v,e) \in V \times E$ and then proceeds by selecting the best enclosing
one. A simple heuristic, which worked well for all our test cases, helps to considerably reduce the number of walks to be processed:
Let
\begin{equation*}
	v_* = \underset{v \in V} \argmin\,d(v) \\
\end{equation*}
and define $W_{v_*}$ as a SEW subject to the constraint
\begin{equation*}
	W_{v_*} \in \tilde{\Pi}_{v_*} = \{ \SP_{v_* ,u} \join uw \join \SP_{w,v_*} : uw \in E \}.
\end{equation*}
Obviously $W_*$ and $W_{v_*}$ do not need to agree in general, as $v_*$ does not have to be traversed by $W_*$.
However, since $v_*$ is the vertex with lowest weight, both walks  differ mainly in a region that has no, or very minor, influence on the total weight and, consequently,
also no significant influence on the condition number.
Actually, $W_*$ and $W_{v_*}$ yielded precisely the same total weight for all functions that we have studied (Fig.~\ref{fig:OneVsAll} compares $W_*$ and $W_{v_*}$ for two typical examples). Using that heuristic,
 the run time of Provan's algorithm improves to $O(|E| + |V| \log|V|)$
because its main part reduces to applying Dijkstra's shortest path algorithm just once. In the case of the grid $\Omega_h$ this bound simplifies to $$O(h^{-2}\log h^{-1}).$$

\subsection{Size of the Grid Domain}
\label{subsec:graphCreation}

The side length $l$ of the square domain supporting $\Omega_h$ has to be chosen large enough to contain a SEW that would approximate an optimal general integration
contour. E.g., if $f$ is entire, we choose $l$ large enough for this square domain to cover the optimal circular contour: $l > 2 r_*$, where $r_*$ is the optimal radius given in \citeasnoun{springerlink:10.1007/s10208-010-9075-z};
a particularly simple choice is $l = 3r_*$.
In other cases we employ a simple search for a suitable value of $l$ by
calculating $W_*$ for increasing values of $l$ until $d(W_*)$ does not decrease substantially anymore. During this search
the grid will be just rescaled, that is, each grid uses a \emph{fixed} number of vertices; this way only the number of search steps enters as an additional factor in the complexity bound. 



\subsection{Multilevel Refinement of the SEW}
\begin{figure}[tbp]
\begin{center}
\begin{minipage}{0.7\textwidth}
\includegraphics[width=\textwidth]{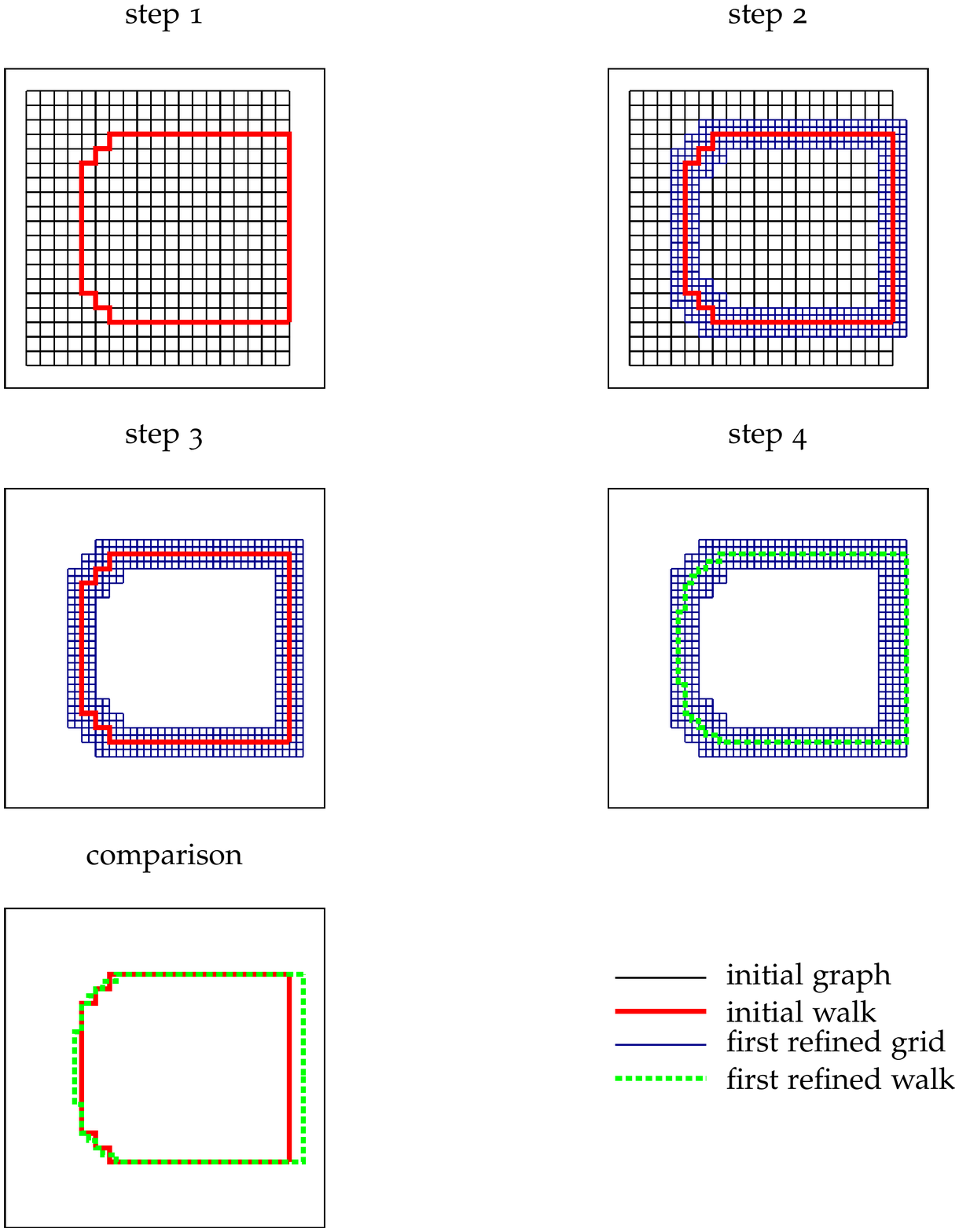}
\end{minipage}	
\end{center}~\\*[-12mm]
\caption{\footnotesize Multilevel refinement of $W_*$ $(f(z)=1/\Gamma(z),\, n=2006)$}
\label{fig:refinement}
\end{figure}

Choosing a proper value of $h$ is not straightforward since we would like
to balance a good approximation of a generally optimal integration contour with a
reasonable amount of computing time. In principle, we
 would construct a sequence of SEWs for smaller and smaller values of $h$
until the total weight of $W_*$ does not substantially decrease anymore. To avoid an undue amount of computational work,
we do not refine the grid everywhere but use an adaptive refinement by
confining it to  a tubular neighborhood of the currently given SEW $W_*$ (see
Fig.~\ref{fig:refinement}):
\begin{itemize}
	\item[1:] calculate $W_*$ within an initial grid;\\*[-3.5mm]
	\item[2:] subdivide each rectangle adjacent to $W_*$ into 4 rectangles;\\*[-3.5mm]
	\item[3:] remove all other rectangles;\\*[-3.5mm]
	\item[4:] calculate $W_*$ in the newly created graph.
\end{itemize}
As long as the total weight of $W_*$ decreases substantially, steps 2 to 4 are repeated.
It is even possible to tweak that process further by not subdividing rectangles that
just contain vertices or edges of $W_*$ having weights below a certain threshold.
By geometric summation, the complexity of the resulting algorithm is 
\[
O(H^{-4} \log{H^{-1}}) + O(h^{-2} \log{h^{-1}})
\]
where $H$ denotes the step size of the coarsest grid and $h = H/2^k$ the step size after
$k$ loops of adaptive refinement. An analogous approach to the constrained $W_{v_*}$-variant of the
SEW algorithm given in §3.2 reduces the complexity further to
\[
O(H^{-2} \log{H^{-1}}) + O(h^{-1} \log{h^{-1}}),
\]
which is close to the best possible bound $O(h^{-1})$ given by the work that would be needed to just list the SEW.

\subsection{Quadrature Rule for the Cauchy Integral}\label{subsect:quadrature} 
\begin{figure}[tbp]
\begin{center}
\begin{minipage}{0.4\textwidth}
{\includegraphics[width=0.99\textwidth]{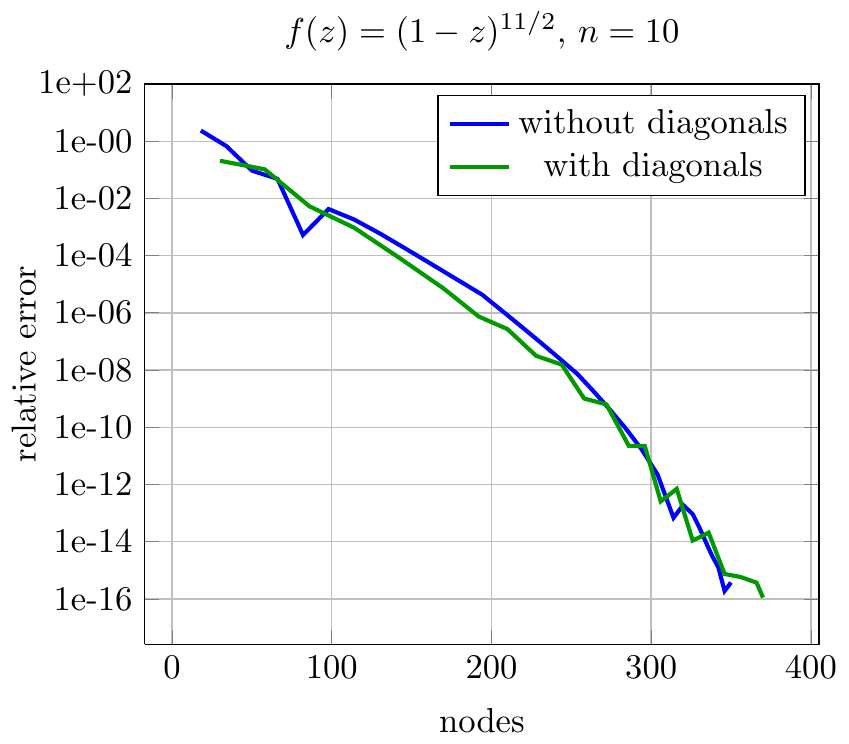}}
\end{minipage}
\hfil
\begin{minipage}{0.4\textwidth}
{\includegraphics[width=0.99\textwidth]{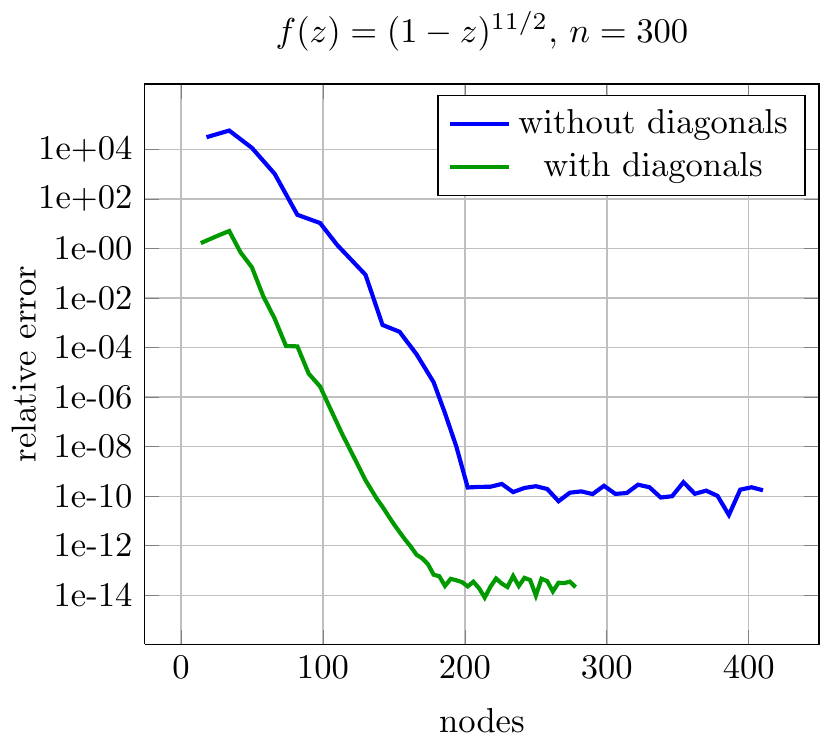}}
\end{minipage}
\end{center}
\vspace*{-3mm}
\caption{\footnotesize Illustration of the spectral accuracy of piecewise Clenshaw--Curtis quadrature on SEW contours for a function
with a branch-cut singularity. For larger $n$, we observe a significant improvement by adding diagonals to the grid. We get to machine precision for $n=10$ and loose about two digits for
$n=300$. (Note that for optimized circular contours the loss would have been about 6 digits for $n=10$ and about 15 digits for $n=300$; cf. Bornemann 2011, Thm.~4.7).}
\label{fig:nodes}
\end{figure}

Finally, after calculation of the SEW $\Gamma = W_*$, the Cauchy
integral (\ref{eq:cauchy}) has to be evaluated by some \emph{accurate} numerical quadrature. We decompose $\Gamma$ into maximally
straight line
segments, each of which can be a collection of many edges. On each of those line segments we
employ Clenshaw--Curtis quadrature in Chebyshev--Lobatto points. Additionally we neglect segments with a weight
smaller than $10^{-24}$ times the maximum weight of an edge of $\Gamma$, since such segments will not contribute to the result within machine precision. This way we not only get \emph{spectral accuracy} but also, in many cases, less nodes as would be needed by the vanilla version of trapezoidal
sums on a circular contour: Fig.~\ref{fig:nodes} shows an example with the order $n=300$ of differentiation but
accurate solutions using just about $200$ nodes which is well below what the sampling condition would require for circular contours
\cite[§2.1]{springerlink:10.1007/s10208-010-9075-z}. Of course, trapezoidal sums would also benefit from some recursive
device that helps to neglect those nodes which do not contribute to the numerical result.

\begin{table}[tbp]
\caption{\footnotesize Condition numbers for some $f(z)$: $r_*$ are the optimal radii given in \protect\citeasnoun{springerlink:10.1007/s10208-010-9075-z};
$W_*$ was calculated in all cases on a $51 \times 51$-grid with $l=3r_*$ (in the last two cases $l$ was
found as in §\protect\ref{subsec:graphCreation}). For $1/\Gamma(z)$, the peculiar order of differentiation $n=2006$ is one of the very rare resonant cases (specific to this entire function) for which circles give exceptionally large condition numbers \protect\citeaffixed[Table~5]{springerlink:10.1007/s10208-010-9075-z}{cf.}. In the
last example, differentiation is for $z=1/\sqrt2$.}\label{tab:1}
\vspace*{0mm}
\centerline{%
\setlength{\extrarowheight}{3pt}
{\footnotesize\begin{tabular}{lrrr}\hline
$f(z)$                                   &   $n$   &      $\kappa(W_*,n)$ & $\kappa(C_{r_*},n)$ \\\hline
$e^z$                                    &   $300$ &      $1.1$           &  $1.0$ \\
$\Ai(z)$                                 &   $300$ &      $1.3$           &  $1.2$ \\
$1/\Gamma(z)$                            &   $300$ &      $1.7$           &  $1.6$ \\
$1/\Gamma(z)$                            &  $2006$ &      $7.8\cdot 10^4$ &  $4.7\cdot 10^4$ \\
\hline
$(1-z)^{11/2}$                           &    $10$ &      $1.4$           &  $5.0\cdot 10^5$ \\
$\exp(1/(1+8z)^{1/5})(1-z)^{11/2}J_0(z)$  &   $100$ &    $7.2 \cdot 10^2$           &  $4.3\cdot 10^{12}$ \\*[0.5mm]\hline
\end{tabular}}}
\label{tab:conditions}
\end{table}

\begin{table}[tbp]
\caption{\footnotesize CPU times for the examples of Table~\ref{tab:1}. Here $t_{W_*}$ and $t_{W_{v_*}}$ denote the
times to compute $W_*$ and $W_{v_*}$ and $t_\text{quad}$ denotes the time to approximate the integral~(\ref{eq:cauchy}) on such a contour by quadrature. (There is no difference between $W_*$ and $W_{v_*}$ from
the point of quadrature,
see Fig.~\ref{fig:OneVsAll}.) In the
last example, differentiation is for $z=1/\sqrt2$. The timings for the grids of size $25\times 25$ and $51\times 51$ match nicely
the $O(h^{-4}\log h^{-1})$ complexity for $W_*$ and the $O(h^{-2}\log h^{-1})$ complexity for $W_{v_*}$. 
}\label{tab:2}
\vspace*{0mm}
\centerline{%
\setlength{\extrarowheight}{3pt}
{\footnotesize\begin{tabular}{lrrrrr}\hline
$f(z)$            & $n$     & grid           & $t_{W_*}$        & $t_{W_{v_*}}$  & $t_{\text{quad}}$  \\ \hline
$e^z$             &  $300$  & $51 \times 51$ &  $4.4 \cdot 10^2$ s & $1.5$ s         & $0.3$ s \\
$\Ai(z)$          &  $300$  & $25 \times 25$ &  $2.1 \cdot 10^1$ s & $0.5$ s         & $1.7$  s \\
$\Ai(z)$          &  $300$  & $51 \times 51$ &  $4.0 \cdot 10^2$ s & $2.1$ s         & $2.1$  s \\
$1/\Gamma(z)$     &  $300$  & $25 \times 25$ &  $2.0 \cdot 10^1$ s & $0.5$ s         & $1.5$  s \\ 
$1/\Gamma(z)$     &  $300$  & $51 \times 51$ &  $3.6 \cdot 10^2$ s & $2.4$ s         & $1.3$  s \\
$1/\Gamma(z)$     & $2006$  & $51 \times 51$ &  $3.6 \cdot 10^2$ s & $2.3$ s         & $3.1$  s \\ \hline
$(1-z)^{11/2}$    &   $10$  & $51 \times 51$ &  $1.4 \cdot 10^3$ s & $5.9$ s         & $0.2$ s \\
$\exp(1/(1+8z)^{1/5})(1-z)^{11/2}J_0(z)$    & $100$    & $51 \times 51$ & $7.0\cdot 10^2$ s & $3.5$ s        & $0.3$ s \\*[0.5mm]\hline
\end{tabular}}}
\label{tab:times}
\end{table}

\begin{figure}[tbp]
\begin{minipage}{0.35\textwidth}
\begin{center}
\includegraphics[width=\textwidth]{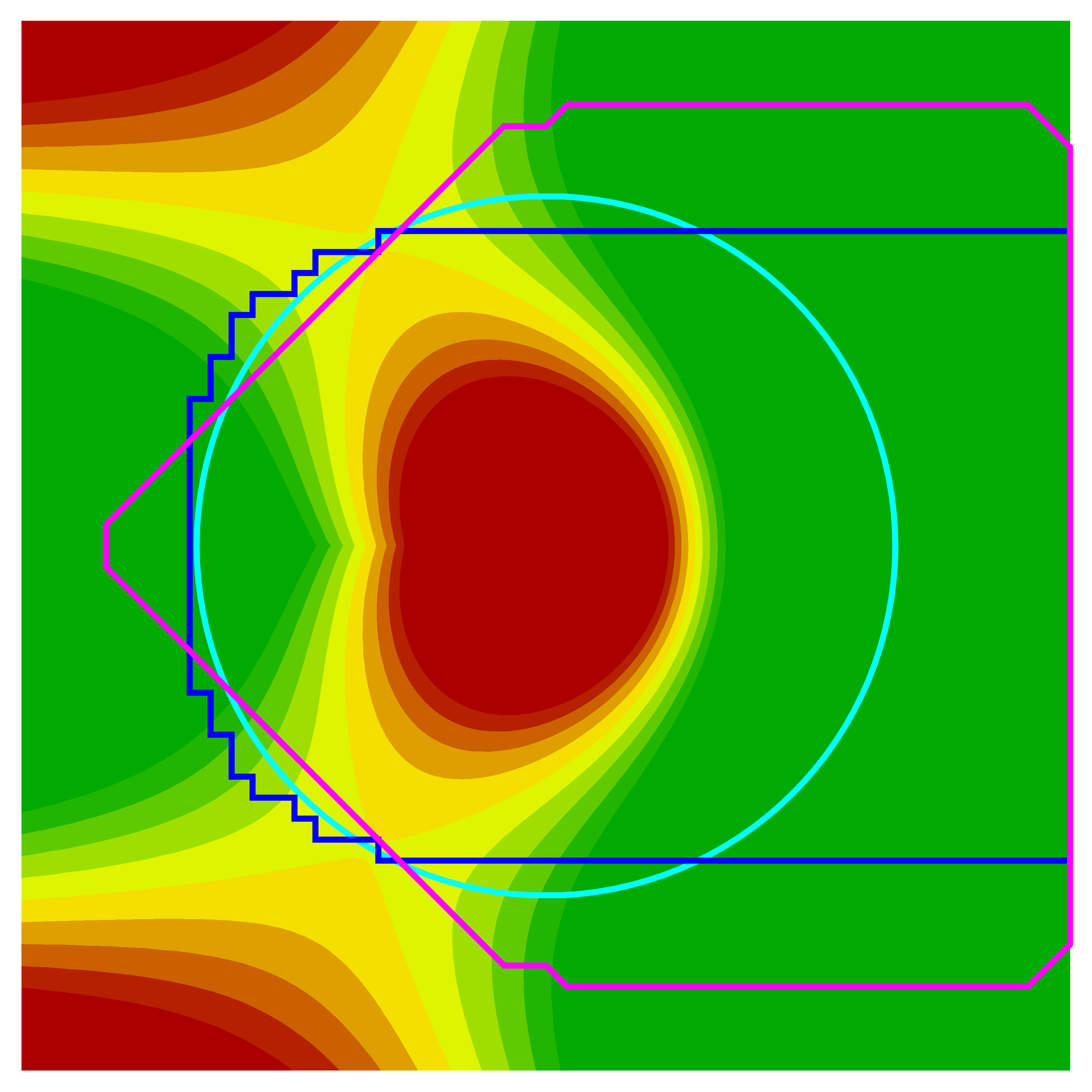}\\*[0.5mm]
{\footnotesize $\Ai(z)$}
\end{center}
\end{minipage}
\hfil
\begin{minipage}{0.35\textwidth}
\begin{center}
\includegraphics[width=\textwidth]{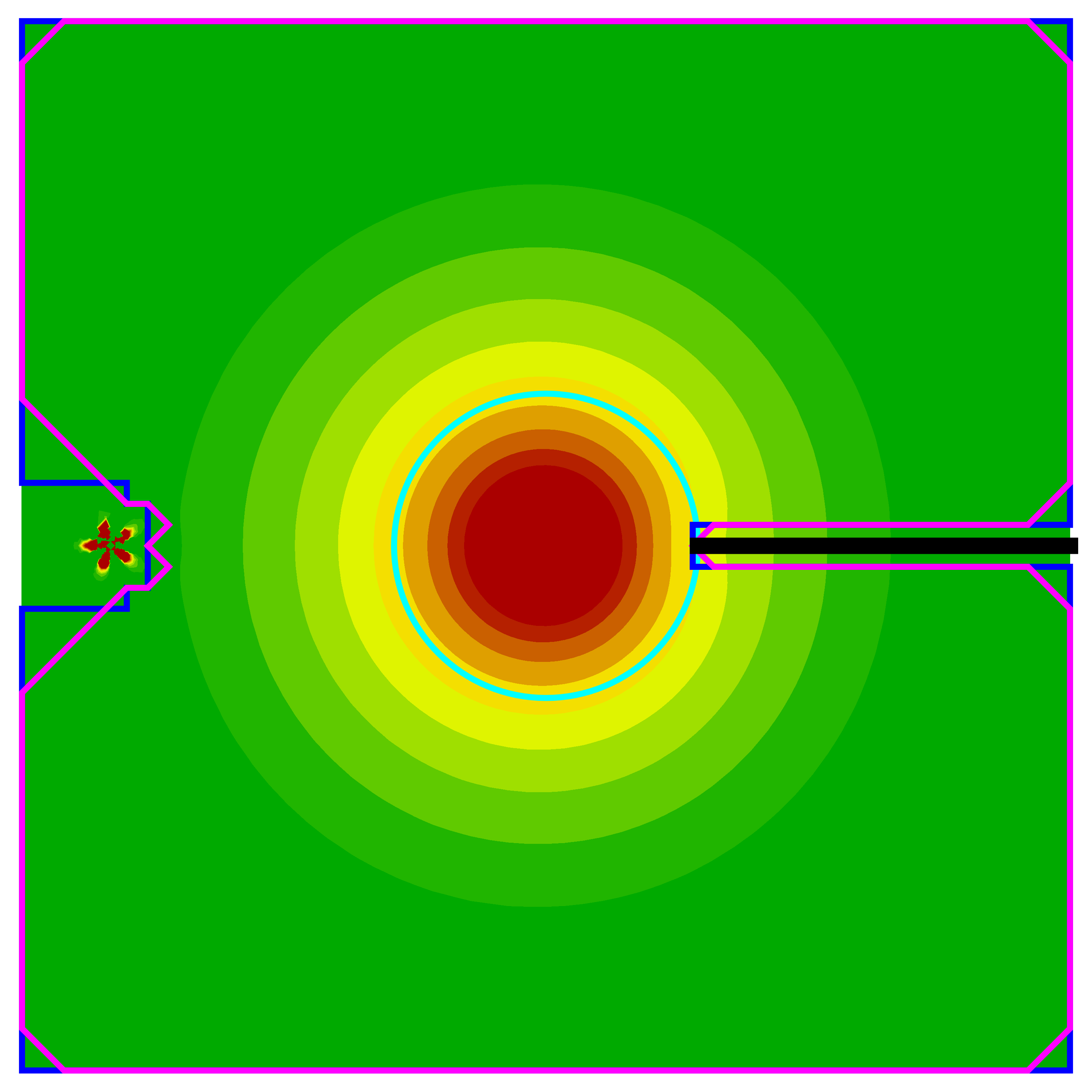}\\*[0.5mm]
{\footnotesize $\exp(1/(1+8z)^{1/5})(1-z)^{11/2}J_0(z)$}
\end{center}
\end{minipage}
	\caption{\footnotesize $W_{v_*}$ (blue: $\Omega_h$ without diagonals, magenta: $\Omega_h$ with diagonals) vs. $C_{r_*}$ (cyan) for some examples of Table~\ref{tab:conditions}:
the color coding shows the size of $\log d(z)$; with red for large values and green for small values. The smallest level shown is the threshold, below of which
the edges of $W_{v_*}$ do not contribute to the first 
significant digits of the total weight.}\label{fig:SewVsCircle}
\end{figure}

\section{Numerical Results}
\label{sec:numer}
Table \ref{tab:conditions} displays condition numbers of SEWs $W_*$  as compared to the optimal circles $C_{r_*}$
for five functions; Table \ref{tab:times} gives the corresponding CPU times
and Fig.~\ref{fig:SewVsCircle} shows some of the contours. ({\em All experiments were done using hardware arithmetic.})
The purpose of these examples is twofold, namely to demonstrate 
that:
\begin{enumerate}
\item the SEW algorithm \emph{matches} the quality of circular contours in cases where the latter are known to
be optimal such as for entire functions;
\item the SEW algorithm is \emph{significantly better} than the circular contours in cases where the latter are known to have
severe difficulties.
\end{enumerate}
Thus, the SEW algorithm is a flexible automatic tool that covers various classes of holomorphic functions in a completely algorithmic fashion; in particular there is 
no deep theory needed to just let the computation run.

In the examples of entire $f$ we observe that $W_*$ and $W_{v_*}$, like the optimal circle~$C_{r_*}$ would do, traverses the saddle points of $d(z)$.
It was shown in \citeasnoun[Thm.~10.1]{springerlink:10.1007/s10208-010-9075-z}  that, for such $f$, the major contribution of the condition number comes
from these saddle points and that circles are (asymptotically, as $n\to \infty$) paths of steepest decent. Since $W_*$ can cross a saddle point
only in a horizontal, vertical, or (if enabled) diagonal direction, somewhat larger condition numbers have to be expected. However, the order of magnitude of the condition number of $C_{r_*}$ is
precisely matched. This match holds in cases where circles give a condition number of approximately $1$, as well as in cases with exceptionally
large condition numbers, such as for $f(z)=1/\Gamma(z)$ in the peculiar case of the order of differentiation $n=2006$ \citeaffixed[§10.4]{springerlink:10.1007/s10208-010-9075-z}{cf.}.

For non-entire $f$, however, optimized circles will be far from optimal in general: \citeasnoun[Thm.~4.7]{springerlink:10.1007/s10208-010-9075-z}
shows that the optimized circle $C_{r_*}$ for functions $f$ from the Hardy space $H^1$ with boundary values in $C^{k,\alpha}$ yields a lower condition number bound of the form
\[
\kappa(C_{r_*},n) \geq c n^{k+\alpha};
\]
for instance, $f(z)=(1-z)^{11/2}$ gives $\kappa(C_{r_*},n) \sim 0.16059\cdot n^{13/2}$.
On the other hand, $W_*$ gives condition numbers that are orders of magnitude better than those of $C_{r_*}$ by automatically following the branch cut at $(1,\infty)$.

\begin{figure}[tbp]
\begin{center}
\begin{minipage}{0.5\textwidth}
{\includegraphics[width=0.99\textwidth]{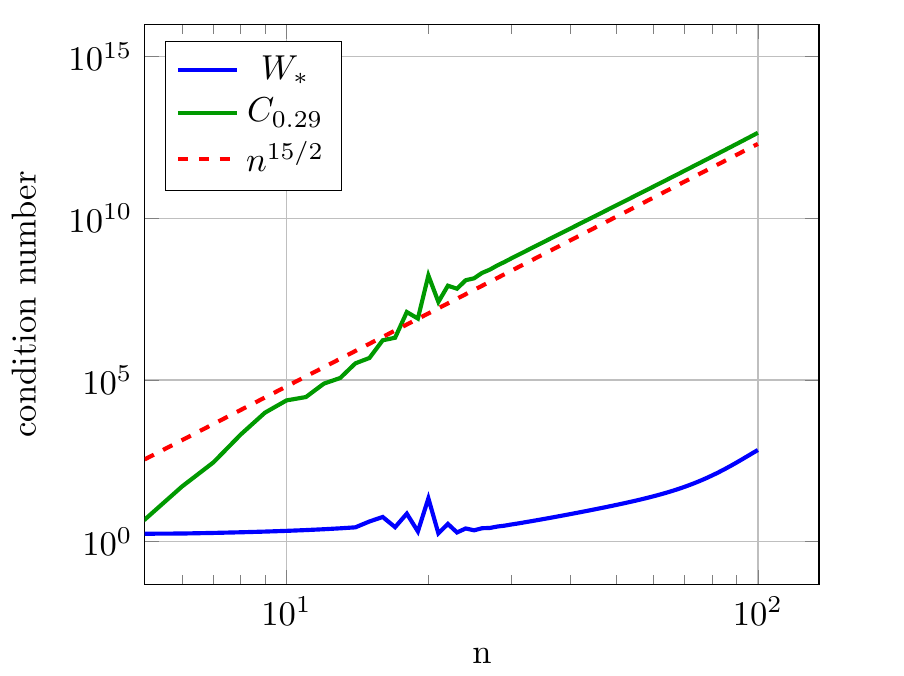}}
\end{minipage}
\end{center}
\vspace*{-3mm}
\caption{\footnotesize An example with essential and algebraic singularities: the condition number of the Cauchy integral  for $\exp(1/(1+8z)^{1/5})(1-z)^{11/2}J_0(z)$ 
for varying order $n$ of differentiation at $z=1/\sqrt 2$; blue: optimal contour $W_*$ in a $51\times 51$ grid graph; green: circular
contour with near optimal radius $r=0.29 \approx 1-1/\sqrt2$; red: prediction of the growth rate from \protect\citeasnoun[Thm.~4.7]{springerlink:10.1007/s10208-010-9075-z}.}
\label{fig:cond}
\end{figure}

The latter example can easily be cooked-up to outperform symbolic differentiation as well: 
using {\em Mathematica 8}, the calculation of the $n$-th derivative
of $f(z)=\exp(1/(1+8z)^{1/5})(1-z)^{11/2}J_0(z)$ at $z=1/\sqrt{2}$ takes already about a minute for
$n=23$ but had to be stopped after \emph{more than a week} for $n=100$. Despite the additional 
difficulty stemming from the combination of an algebraic and an essential singularity at $z=-1$, the $W_{v_*}$ version of the SEW
calculates this $n=100$ derivative to 
an accuracy of 13 digits in less than $4\,$s; whereas optimized circular contours would give only about 3 correct digits here (see Fig.~\ref{fig:cond}).

While many more such numerical experiments would demonstrate that reasonably small condition numbers are obtainable in general,\footnote{The software is provided as a supplement to the e-print version of this paper: \href{http://arxiv.org/e-print/1107.0498v3}{\tt arXiv:1107.0498}.}
the study of rigorous condition number bounds for the SEW has to be postponed to future work. 

\bibliographystyle{kluwer}
\bibliography{cauchyint}
\end{document}